# A Survey of Hamilton's Program for the Ricci Flow on 3-manifolds

Bennett Chow*

Department of Mathematics
University of California, San Diego
La Jolla, CA 92093

## 1 Introduction

This is a purely expository article on Riemannian Ricci flow with emphasis on dimension three. None of the results in this paper are due to the author. In view of Hamilton's vast works, Yau's influence on the field, and Perelman's recent inspiring paper on Ricci flow, there may be more interest in the area. The author has only partial knowledge of this field but hopefully this paper will help nonspecialists and graduate students to follow some of the recent developments. Unfortunately, this paper does not include the very recent work of Perelman [P], which hopefully will be addressed at some later time.[1] A previous expository article on Ricci flow was written by H.-D. Cao and the author [CaCh], we refer to the bibliography in there for some more references. In [H7], Hamilton included a survey and then proved a plethora of new results on Ricci flow with particular emphasis on singularity analysis.

## 2 Ricci flow and geometrization

Let $M^3$ be a closed, orientable 3-manifold. Thurston's Geometrization Conjecture says that there exists a finite collection of disjoint, embedded 2-spheres and 2-tori such that after cutting $M^3$ along these surfaces and capping the boundary 2-spheres by 3-balls, the interior of each component of the resulting 3-manifold admits a complete locally homogeneous metric [T1]. The three most well-known homogeneous geometries are the constant sectional curvature $+1, 0$ and $-1$ geometries, corresponding to the simply-connected models: sphere $S^3$, euclidean space $\mathbb{R}^3$, and hyperbolic space $\mathbb{H}^3$, respectively. These geometries are isotropic, that is, they look the same in every direction. In particular, the isometry groups of these spaces act transitively on their orthonormal frame bundles. After that, there are the product geometries with models: $S^2 \times \mathbb{R}$ and $\mathbb{H}^2 \times \mathbb{R}$ ($\mathbb{R}^2 \times \mathbb{R}$ is the same as $\mathbb{R}^3$,) which are the products of surface geometries and the line. Finally, the least isotropic geometries are modelled on 3-dimensional unimodular Lie groups: Heisenberg (also called Nil,) Solv, and $\widetilde{SL_2\mathbb{R}}$, which may be considered as twisted products of $\mathbb{R}$ with $\mathbb{R}^2$ and $\mathbb{H}^2$.

---

*Partially supported by NSF grant DMS 0203926.
[1] Due to the short amount of time the author has had to prepare this article, it is very incomplete and most likely inaccurate at places. We apologize for any errors and omissions. We have also not included the figures, due to a difficulty in uploading .wmf files.

Now how does Hamilton's Ricci flow fit in with all of this? His idea is to start with a closed, orientable, differentiable 3-manifold $M^3$, endow it with a smooth Riemannian metric (this is always possible,) and then run the volume normalized Ricci flow.[2] The philosophy is that from the analytic viewpoint, the Ricci flow is the heat equation for Riemannian metrics, and hence will smooth out the metric as much as possible. Furthermore, from the geometric viewpoint, the Ricci flow is very natural: it evolves the metric in the direction of an intrinsically defined curvature tensor. Because of this, the Ricci flow is invariant under the diffeomorphism group of the manifold $M^n$: if $g(t)$ is a solution to the Ricci flow, then for any fixed diffeomorphism $\phi : M^n \to M^n$, the family $\phi^* g(t)$ is also a solution to the Ricci flow. Hence any symmetries of the initial metric are preserved under the flow. In particular, if the initial metric is homogeneous, then the solution metrics are homogenous for all positive time. In general, one expects that the Ricci flow tries to make the metric more (at least not less) homogeneous and isotropic (see [IJ] for a study of the Ricci flow of homogenous metrics).

If the volume normalized Ricci flow converges, the limit metric is a fixed point of the flow and hence a constant sectional curvature metric. Thus $M^3$ admits a constant sectional curvature geometry and is diffeomorphic to a space form. Spherical (curvature $+1$) space forms and flat manifolds (Bieberbach) are classified in J. Wolf's book [W]. Hyperbolic manifolds (curvature $-1$) are well-understood by the work of Thurston [T2], [T3], [T4] and others.

If the curvature of the solution metric stays uniformly bounded (such solutions are called nonsingular,) while the flow does not converge either, then we still expect something nice to happen due to the uniform bound on curvature. In fact, Hamilton has proved that in this case $M^3$ admits a geometric decomposition [H10]. What happens is that as time approaches infinity, if the metrics do not collapse in the sense of Cheeger-Gromov (see [CG1] and [CG2]), then disjoint pieces of $M^3$ limit to a finite collection of complete, noncompact, finite volume, hyperbolic manifolds. These hyperbolic limits, after truncating the ends, can be embedded in $M^3$ so that their boundaries (consisting of 2-tori) are incompressible in $M^3$ and are joined in $M^3$ by a 3-dimensional submanifold with 2-tori boundary whose components admit Cheeger-Gromov collapsing metrics (since collapsing manifolds are graph manifolds, they also admit geometric decompositions.) Geometrically, for $t$ large enough, there exists a finite number of disjoint, embedded, incompressible 2-tori $\cup T_i(t)$ depending smoothly on $t$ which decomposes $M^3$ into pieces, each of which either approaches a complete, noncompact, finite volume, hyperbolic 3-manifold or collapses as $t$ tends to infinity.

---

[2] Given smooth initial data, a solution always exists for short time [H1]. See [D] for a simplified proof of this.

When singularities do occur (i.e., the curvature becomes unbounded,) they are expected to be topologically and geometrically simple. Somewhere near the singularity the 3-manifold is believed to look like the shrinking round product cylinder $S^2 \times \mathbb{R}$ (or one of its quotients.) This will allow for topological-geometric surgeries to be performed right before the singularity time. Starting with the new Riemannian 3-manifold, we start the Ricci flow again and repeat as many times as necessary the process of performing topological-geometric surgeries right before any singularity times. In the following sections we shall give an intuitive description of when, where, and how to perform these surgeries. Roughly speaking, surgeries can be performed on all necks in the 3-manifold which are geometrically close enough to the quotient of a round product cylinder. These surgeries are performed at the first times they reduce the maximum absolute curvature by a large enough factor.

Now Hamilton's deep result is that if we can get a good enough classification of the possible limits of dilations about singularities[3], then we can perform a geometric-topological surgery on the singular solution $(M^3, g(t))$ to the volume normalized Ricci flow, where we immediately throw away resulting components which are topologically $S^3/\Gamma$ or $(S^2 \times \mathbb{R})/\Gamma$, at some time $\tau_1 < T$, which depending on the choice of $\tau_1$, decreases the maximum absolute curvature of $(M^3, g(\tau_1))$ by an arbitrarily large factor. His even deeper theory is that such a classification is hoped to imply (by using a contradiction argument) that there exist a finite number of such surgeries after which the final solution to the volume normalized Ricci flow is nonsingular (on a closed 3-manifold which may be empty.) Thus, such a classification of limits of dilations about singularities would imply the geometrization all closed 3-manifolds. In the following sections we describe some of the intuition behind the results above.

## 3 Neck pinches

### 3.1 Intuition

The intuition and hope we have about the possible types of singularities of the Ricci flow which may form on a closed, orientable 3-manifold $M^3$ is that they either look like a neck pinching off or a degenerate neck pinching off. Since degenerate neck pinches are expected to be non-generic, for simplicity we now restrict ourselves to neck pinches. First of all, a neck is an open subset $N$ of $M^3$ which topologically is diffeomorphic to a quotient of $S^2 \times \mathbb{R}$. If $N$ is also closed (i.e., compact,) then $N = M^3$, and $M^3$ is classified as an orientable, compact quotient of $S^2 \times \mathbb{R}$. These are $S^2 \times S^1$, $\mathbb{R}P^3 \# \mathbb{R}P^3$, and $S^2 \tilde{\times} S^1$: the topologically unique nontrivial $S^2$ bundle over $S^1$. Note that all of these spaces admit geometric structures. If $N$ is noncompact, then $N$ is diffeomorphic to either $S^2 \times \mathbb{R}$ or $\mathbb{R} \tilde{\times} \mathbb{R}P^2$: the unique nontrivial $\mathbb{R}$ bundle over $\mathbb{R}P^2$. Note also that $\mathbb{R} \tilde{\times} \mathbb{R}P^2 \cong \mathbb{R}P^3 - \{\text{ball}\}$.

---

[3]According to Hamilton, the fundamental estimate missing to obtain this classification is the Harnack inequality for the Ricci flow for any initial metric on a closed 3-manifold.

Roughly speaking, given a solution to the Ricci flow $\left(M^3, g\left(t\right)\right)$ on $[0, T)$, a neck pinch in $M^3$ at the singularity time $T$ is a time-dependent neck $N\left(t\right)$ with the induced metric $g\left(t\right)$ which asymptotically approaches the quotient of a shrinking round product cylinder $S^2 \times \mathbb{R}$ as $t \nearrow T$. Thus we see a thin and relatively long cylinder forming.

Now this of course is not the complete picture of the forming singularity, which would require a detailed analysis of the asymptotics. On the other hand, we are interested only in understanding enough about the singularity so that we can perform surgeries effectively. Intuitively, we expect that the forming neck actually looks like a very thin and long cylinder which opens up very slowly (so that it looks like horn(s) at the end(s).)

That is, the neck is not geometrically close to just one cylinder but close to many cylinders of varying radii depending on where on the neck you are: a point near the center of the neck is contained in a subneck of a thin radius, and as the point moves away from the center and toward the end(s), it is contained in subnecks of larger and larger radii. In particular, we expect the following somewhat more precise formulation. As we approach the singularity, the neck *conformally* approaches a round product cylinder. That is, there exists positive functions $\psi\left(t\right)$ on the necks $N\left(t\right)$ such that the metric $\psi\left(t\right) g\left(t\right)$ approaches a noncompact quotient of the round product cylinder. The functions are approximately constant on the $S^2$ slices and very slowly increase as we go out from the center.

## 3.2  The dilation argument

Now why do we expect all of this to be true? The basic reason is as follows. If a neck pinches, then there exists a sequence of points and times $(x_i, t_i)$ with $x_i \in M^3$ and $t_i \nearrow T$ such that if for each $i$ we center microscopes at $(x_i, t_i)$ and magnify the metric so that its absolute curvature at $(x_i, t_i)$ becomes 1, then the dilated metrics converge to the quotient of a round product cylinder. Now we assume the absolute curvatures at $(x_i, t_i)$ are comparable to the maximum absolute curvatures of $g(t_i)$ on $M^3$ (this should be true for a neck pinch but not for a degenerate neck pinch.) That is, about the sequence of points and points $(x_i, t_i)$ with absolute curvatures comparable to the spatial maximums, asymptotically we see a round product cylinder. Now the point is from Hamilton's theory, it follows that if we magnify about *any* sequence of points and times $(y_i, u_i)$ with $y_i \in M^3$ and $u_i \nearrow T$ with absolute curvatures comparable to the spatial maximums, we should also see a round product cylinder. Note that we cannot obtain a compact quotient of a cylinder as a limit by the resolution of Ilmanen and Knopf of a conjecture of Hamilton [IK]. See also [GIK] for the application of center manifold analysis to the stability of the Ricci flow in some special cases.

Choosing any points with absolute curvatures varying from the spatial maximum to any comparable value, as long as the time is close enough to the singularity time $T$, we see approximate round product cylinders. As a point moves from where the absolute curvature is maximum to any point where the curvature is comparable to the maximum, the corresponding approximate round product cylinders should overlap and connect up to form a very long cylinder which is conformally an approximate round product with the radius opening up as the curvature of the points decrease. This is why near the singularity time, the neck should look like one of the pictures in diagram 4. In fact the above reasoning should actually imply the following. Given any positive numbers $c$ and $\varepsilon$ no matter how small, there should exist a time $t_c < T$ such that any point $(x, t)$ with $t \geq t_c$ and $|Rm|(x, t) \geq c \max_{M \times \{t\}} |Rm|$ is the center of a long neck which is $\varepsilon$-close to a round product cylinder (we'll be more precise about this later.) That is, for any $c > 0$ no matter how small, for times close enough to the singularity time, we see approximate round product cylinders of every scale $c$-comparable to the spatial maximum absolute curvature. This means that if we dilate the center radius of the neck to be 1, near the singularity time, the conformal approximate round product cylinder opens to an arbitrary large radius ($\sim 1/\sqrt{c}$) and the curvature at the end(s) of the cylinder are arbitrarily small compared to the curvature at the center ($\sim c$.)

# 4 Surgery

## 4.1 How to perform it

Now that we have a good intuitive picture of a neck pinch, how should we perform the topological-geometric surgery? The answer is not difficult now that we know so much about what to expect. First consider a neck pinch of the topological type $S^2 \times \mathbb{R}$ (as we shall see, the $\mathbb{R} \tilde{\times} \mathbb{R}P^2$ case is analogous.) Given a neck $N(t)$ at some time $t$ near $T$, which we identify with $S^2 \times [-L, L]$, and a positive number $\delta$, topologically the surgery is performed by removing the subcylinder $S^2 \times [-L + \delta, L - \delta]$ and capping off each of the two 2-spheres by a 3-ball. If the new 3-manifold $\hat{M}^3$ is not connected, then $\hat{M}^3$ has two connected components $M_1$ and $M_2$ and $M^3 \cong M_1 \# M_2$. If $\hat{M}^3$ is connected, then $M^3 \cong \hat{M}^3 \# (S^2 \times S^1)$.

Geometrically the new metric $\hat{g}(t)$ on $\hat{M}^3$ is obtained from $g(t)$ on $M^3$ by taking a spherically symmetric positive sectional curvature metric on the 3-ball such as the standard metric on a hemisphere in $S^3$ and blending this in with $g(t)$ on each of the two ends $S^2 \times [-L, -L + \delta]$ and $S^2 \times [L - \delta, L]$. Before we can do this blending, we need to multiply the metrics $g(t)$ on the ends by conformal factors $\xi(t)$ which are constant on the slices $S^2 \times \{s\}$ and which make the metrics $\xi(t) g(t)$ positively curved and such that the 2-spheres $S^2 \times \{s\}$ shrink as $s$ moves toward the inner endpoint ($-L + \delta$ or $L - \delta$.) This allows for the blending to be done in a way in which the part of the 3-manifold where the surgery takes place has positive sectional curvature.

On the other hand, if the neck pinch is of the topological type $\mathbb{R} \tilde{\times} \mathbb{R}P^2$, then we identify the neck $N(t)$ with $[-L, L] \tilde{\times} \mathbb{R}P^2 \cong [0, L] \times S^2/ \sim$, where $(0, x) \sim (0, -x)$. We perform the surgery topologically by capping off the end $S^2 \times [L-\delta, L]$ with a 3-ball and geometrically by blending the metrics in a way exactly analogous to the previous case. Throwing away the piece $[0, L-\delta] \times S^2/ \sim$ yields a new 3-manifold, which we call $M^*$. We then have $M^3 \cong M^* \# \mathbb{R}P^3$ since $[0, L-\delta] \times S^2/ \sim$ is diffeomorphic to $\mathbb{R}P^3 - \{3\text{-ball}\}$.

## 4.2 Dropping the maximum curvature and obtaining non-singular solutions

For the ideas in this section it is important to recall that whenever we perform a surgery we immediately throw away resulting components which are topologically $S^3/\Gamma$ or $(S^2 \times \mathbb{R})/\Gamma$.

As before, we claim that given any positive constant $\Gamma$ no matter how large, we expect to be able to perform one of the above surgeries at some time $t$ a little before $T$ such that the maximum absolute curvature of the metric after surgery decreases by the factor $\Gamma$. The reason is that given $\Gamma$, we just choose the constant $c$ in section 3.2 on the order of $1/\Gamma$, and perform the surgery at time $t_c$ or after. The surgery is performed near the end(s) of the conformal approximate round product cylinder(s) where the absolute curvature is approximately $c$ times the maximum absolute curvature on $M^3$. Hence, regions where the absolute curvatures are at least $c$ times their maximums are removed. At a somewhat more detailed level, we present the intuition for why there should exist a finite number of surgeries, after which the solution is nonsingular. A simpler version of this argument provides some of the details the above claim that surgery should reduce the maximum absolute curvature by an arbitrary large factor.

We claim that there should exist a certain set of parameters (previously we oversimplified these parameters to just one small positive number $\varepsilon$) measuring how close to a round product cylinder we require necks to be and by how large a factor ($\Gamma$) we want the standard surgery on maximal necks of such type to reduce the maximum absolute curvature, such that the sequence of first times $\{\tau_j\}$ when such surgeries can be performed is finite and that the final solution is nonsingular (implicit in this statement is that for a singular solution, such a surgery exists in the first place.) If such a set of parameters does not exist, then no matter how good we choose these parameters, either the sequence of first times $\{\tau_j\}$ when such surgeries can be performed is infinite or $\{\tau_j\}$ is finite and the final solution is singular. The second case is equivalent to saying that for the final singular solution, it is not possible to perform a surgery which reduces the curvature by the given factor. Since the argument in the first case essentially contains the argument in the second case, we just consider the first case where $\{\tau_j\}$ is infinite.

Let $\tau_0 \doteq 0$ and $(M_j, g_j(t))$, $t \in [\tau_j, \tau_{j+1}]$, be the $j$th solution to the normalized Ricci flow where $(M_{j+1}, g_{j+1}(\tau_{j+1}))$ is obtained from $(M_j, g_j(\tau_{j+1}))$ by surgery and the zeroth solution is the original solution. An important fact to keep in mind is that the curvatures of the solutions cannot stay uniformly bounded since then each surgery will decrease the volume by at least a fixed amount, but we started with a closed 3-manifold so that the initial volume is finite. Now the idea is to find a new infinite sequence of points and times $\{(p_j, \sigma_j)\}$ and take a limit of dilations about this sequence and somehow obtain a contradiction. First we address how to choose this sequence of points and times. We choose $\sigma_j \in (\tau_j, \tau_{j+1})$ to be the first time such that the maximum absolute curvature of $g_j(\sigma_j)$ is equal to $g_{j-1}(\tau_j)$ (recall that the surgery drops the maximum by a large factor and also that the curvature of the collection of solutions is unbounded.) The point is that since $\sigma_j < \tau_{j+1}$, at the time $\sigma_j$ the standard surgery on all necks in $(M_j, g_j(\sigma_j))$ does not decrease the maximum absolute curvature by the preset large factor. This means that there is a sequence of points $p_j \in M_j$ with absolute curvatures uniformly comparable to the maximums at the times $\sigma_j$ such that the points $p_j$ are not contained in the necks at the times $\sigma_j$ nor in the regions of $M_j$ which the standard surgeries remove.

Now take the limit of dilations about a subsequence of $\{(p_j, \sigma_j)\}$ (making the big assumption that we have an injectivity radius estimate.) There are three cases to consider. In the first case, the limit $M_\infty$ is a spherical space form. This implies that one of the components of $M_j$ is diffeomorphic to $M_\infty$ for $j$ large enough, contradicting the fact that we immediately throw away any spherical space forms components when they arise. The second case is where $M_\infty$ is a quotient of a round product cylinder. If the quotient is compact, then components of $M_j$ are diffeomorphic to $M_\infty$ for $j$ large enough, contradicting the fact that we also immediately throw away any $\left(S^2 \times \mathbb{R}\right)/\Gamma$ components. If the quotient is noncompact, then for $j$ large enough the points $(p_j, \sigma_j)$ are in the center of regions in $M_j$ arbitrarily close to noncompact quotients of round product cylinders. This contradicts the assumption that the points $(p_j, \sigma_j)$ are not removed by the standard surgeries for all $j$.

The third case is where the limit $(M_\infty, g_\infty(t))$ is noncompact with strictly positive sectional curvature and close to a self-similar solution. In this case, one expects that the limit is close to the quotient of a round product cylinder at spatial infinity (what we picture is similar to the limit of dilations about the tip of a degenerate neck pinch, which we discuss in the next installment.) In any case, we expect that $M_\infty$ is topologically $\mathbb{R}^3$ and geometrically regions near infinity in its end are arbitrarily close to round product cylinders (i.e., the quotient is trivial.)

In the solutions $(M_j, g_j(t))$, $t \in [\tau_j, \tau_{j+1}]$, these regions will correspond to necks which surgery will remove. However, the points $p_j \in M_\infty$ correspond to points in the part of the manifolds $M_j$ which become diffeomorphic to $S^3$ after surgery (since $M_\infty \cong \mathbb{R}^3$.) This means that the standard surgeries at the times $\sigma_j$ remove $p_j$ for $j$ large enough, which is a contradiction. Thus we expect that after a finite number of surgeries, a nonsingular solution arises.

This completes a loose description of Hamilton's program for the Ricci flow on closed 3-manifolds as an approach to Thurston's Geometrization Conjecture.

# 5 Singularity analysis

Now we begin a slightly more detailed description of the arguments behind the intuitive picture of what we expect in dimension three which we gave in the previous sections. In general, it is difficult to analyze and classify the singularities of solutions to nonlinear partial differential equations. The situation for the Ricci flow equation is no exception, especially since the equation is a system and the solutions are metrics, which means the ambient space does not have a rigid structure such as is the case for euclidean space or when a background metric is fixed on a manifold. As we remarked earlier, the Ricci flow is natural geometrically since the Ricci tensor is intrinsic and the Ricci flow equation is invariant under the diffeomorphism group of the manifold. From the analytic viewpoint, the Ricci flow is very similar to the heat equation (more than one would expect from it just being a second order parabolic quasilinear system!) Just about any geometrically defined scalar or tensor (such as curvature, the length of a geodesic, and the area of a minimal surface) satisfies a nonlinear heat-type equation. This enables one to use a powerful analytic tool: the maximum principle. The application of the maximum principle is behind most of the understanding of singularities of the Ricci flow.

# 6 Singularity types

Given a solution $(M^n, g(t))$ to the Ricci flow on a closed manifold $M^n$ and time interval $[0, T)$, we say that a singularity forms at time $T$ if

$$\sup_{M \times [0,T)} |Rm| = \infty.$$

Singularities are categorized into two distinct types:, depending on the rate of blow up of the curvature:

Type I. $\sup_{M \times [0,T)} (T-t) |Rm| < \infty.$

Type II. $\sup_{M \times [0,T)} (T-t) |Rm| = \infty.$

The prime example of a Type I singularity is that of a round sphere shrinking to a point. If $g(0)$ is a constant positive sectional curvature metric on $S^n$, then the solution to the Ricci flow is $g(t) = \left(1 - \frac{2}{n} R(0) t\right) g(0)$ and $T = \frac{n}{2} R(0)^{-1}$. This expression for the solution metric also holds when $g(0)$ is an Einstein metric with positive scalar curvature. Note also that if we start with a product manifold $(M^n \times N^m, g(0) \times h(0))$ where $g(0)$ is an Einstein metric with positive scalar curvature on $M^n$ and $h(0)$ is a Ricci flat metric on $N^m$, then the solution to the Ricci flow is $g(t) \times h(0)$. The solution metric is the product of a shrinking positive Einstein metric with a Ricci flat metric. The singularity is again Type I. In dimension three, the above examples of Type I singularities correspond to shrinking spherical space forms $\left(S^3/\Gamma, g(t)\right)$ and quotients of shrinking cylinders $\left(\left(S^2 \times \mathbb{R}\right)/\Gamma, (g(t) \times h(0))/\Gamma\right)$. The importance of these two examples is that they are homothetically shrinking solutions, that is, the metric at time $t$ is a constant (depending on $t$) multiple of the initial metric, where the constant decreases to zero as $t$ approaches the final time. As we shall see later (section 10,) most, if not perhaps all, dilations about a Type I singularity on a closed 3-manifold tend to these special solutions.

Type II singularities are harder to picture. As we shall see later, these are expected to correspond to so-called degenerate neck pinches. The corresponding special solutions which appear as limits of sequences of dilations about Type II singularities are solutions which are time independent (modulo the action of the diffeomorphism group.) Besides the Ricci flat solutions, which are stationary under the Ricci flow, these solutions are on noncompact manifolds. Two important examples on noncompact 3-manifolds are $\mathbb{R}^3$ with a certain rotationally symmetric metric with positive sectional curvature, and $\mathbb{R}^3 = \mathbb{R}^2 \times \mathbb{R}^1$ with the product metric of a 'cigar soliton' solution and the line. The cigar soliton solution is a time independent solution to the Ricci flow on $\mathbb{R}^2$ which is given by the initial metric $g\left(0\right) = \left(dx^2 + dy^2\right)/\left(1 + x^2 + y^2\right)$. It has positive curvature and is asymptotic to a cylinder of radius 1.

# 7  Dilations

Before we begin our analysis, we first recall (in somewhat more detail than in the previous article) how to take dilations of the solution about a singularity. Choose a sequence of points and times $(x_i, t_i) \in M \times [0, T)$ with $t_i \nearrow T$ such that $K_i \doteq |Rm|\left(x_i, t_i\right) \geq c \cdot \max_{M \times [0, t_i]} |Rm|$ for some constant $c > 0$ independent of $i$. This condition guarantees that when we translate time $t_i$ to time 0 and rescale the metrics and time so that the new solutions $\tilde{g}_i\left(t\right) \doteq K_i \cdot g\left(t_i + K_i^{-1}t\right)$ satisfy $\left|\widetilde{Rm_i}\right|\left(x_i, 0\right) = 1$, the curvatures of the new solutions are uniformly bounded (by $1/c$) at all times before 0. By the Gromov type compactness theorem for solutions of the Ricci flow [H6],[4] provided we have a uniform injectivity radius estimate for the metrics $\tilde{g}_i\left(0\right)$ at the points $x_i$, there exists a subsequence such that the pointed solutions $\left(M^3, \tilde{g}_i\left(t\right), x_i\right)$ converge uniformly in $C^\infty$ on compact sets to a complete solution $\left(M^3_\infty, g_\infty\left(t\right), x_\infty\right)$ to the Ricci flow. Note that this limit solution is defined on an interval $(-\infty, \omega)$, where $\omega > 0$ (such solutions are called *ancient*.) When we speak of classifying singularities, what we really mean is classifying the possible limits of such dilations (and certain other limits of sequences of dilations where the absolute curvature at $(x_i, t_i)$ is not uniformly comparable to its maximum at time $t_i$.) One of the main obstacles is obtaining the necessary injectivity radius estimate so that we can take such limits (more on the injectivity radius estimate in a later article.)

# 8  Type I

Type I singularities of solutions to the Ricci flow on closed, orientable 3-manifolds are now well-understood. First, Hamilton has obtained an isoperimetric esti-mate for solutions forming Type I singularities. This estimate implies the necessary injectivity radius estimate. In particular, there exists a constant $c > 0$ such that for any point $(x, t)$ we have

$$\operatorname{inj}(x, t) \geq c\sqrt{T - t}.$$

---

[4]See [Lu] for an extension to orbifolds and [Gl] for an extension to collapsing sequences.

Since for a Type I singularity we have $c_1 / (T - t) \leq \max_{M \times \{t\}} |Rm| \leq c_2 / (T - t)$ for some constants $c_1 > 0$ and $c_2 < \infty$, this is equivalent to

$$\max_{M \times \{t\}} \sqrt{|Rm|} \cdot \operatorname{inj}(x, t) \geq c > 0.$$

This yields a uniform injectivity radius estimate for the sequence of dilated metrics $\tilde{g}_i(0)$ defined above. Hence we can apply the compactness theorem to obtain the desired limit solution $(M_\infty^3, g_\infty(t), x_\infty)$. By a clever maximum principle argument, Hamilton has shown that for any singularity there exists a sequence of points and times and a corresponding limit solution which is either a shrinking spherical space form $(S^3/\Gamma, g(t))$ or a quotient of a shrinking cylinder $((S^2 \times \mathbb{R})/\Gamma, (g(t) \times h(0))/\Gamma)$, which are the two model Type I singularities we discussed in section 6.

In the first case, the limit $M_\infty^3 \cong S^3/\Gamma$ is compact. This implies that the original manifold $M^3$ is diffeomorphic to $S^3/\Gamma$ (since $M^3 \cong M_\infty^3$.) In this case we have classified $M^3$ as a quotient of $S^3$ by a finite subgroup of $O(4)$ acting freely and properly discontinuously.

In the second case, the limit manifold $M_\infty^3 \cong (S^2 \times \mathbb{R})/\Gamma$ may be either compact or noncompact. If $M_\infty^3$ is compact, then again $M^3 \cong M_\infty^3$ and we have classified $M^3$ as either $S^2 \times S^1$, $\mathbb{R}P^3 \# \mathbb{R}P^3$, or $S^2 \tilde{\times} S^1$ (the topologically unique nontrivial $S^2$ bundle over $S^1$.)

If $M_\infty^3$ is noncompact, then it is diffeomorphic to either $S^2 \times \mathbb{R}$ or $\mathbb{R} \tilde{\times} \mathbb{R}P^2$ (the unique nontrivial $\mathbb{R}$ bundle over $\mathbb{R}P^2$.) Geometrically, this corresponds to the formation of a neck in $M^3$ at time $T$.

We conclude that if a Type I singularity forms, then $M^3$ is either diffeomorphic to a quotient of either $S^3$ or $S^2 \times \mathbb{R}$, or a neck forms at the singularity time (i.e., there exists a noncompact cylindrical limit.)

The following question appears to be still open. Can one classify all Type I limits? Perhaps we may conjecture that the above limits are the only possibilities (see [CL]).

We also remark that recently Angenent and Knopf [AK] have studied the existence of neck pinches for compact rotationally symmetric solutions in part based on previous work of Miles Simon in the noncompact case [S].

# 9   Type II

At the moment, there is only a partial understanding of Type II singularities in dimension three. In particular, one expects that the injectivity radius estimate still holds and that one of the singularity models, namely the cigar soliton does not to occur as a limit. The validity of both of these statements is unknown. However, if both of these statements are true, then for any Type II singularity, there exists a sequence of points and times $(x_i, t_i) \in M \times [0, T)$ with $t_i \nearrow T$ such that the limit of a subsequence of dilations is a round product cylinder $\left(M_\infty^3, g_\infty(t), x_\infty\right) = \left(S^2 \times \mathbb{R}, g(t) \times h(0)\right)$. This detects the neck that we wanted to find. We shall outline Hamilton's proof of this in the following sections.

Now for a Type II singularity, we expect that $\sup_{M \times \{t_i\}} (T - t_i) |Rm| \to \infty$ for any sequence $t_i \nearrow T$. On the other hand, if the limit is a quotient of a cylinder which is an ancient Type I solution, then we expect that the sequence of points and times $(x_i, t_i)$ satisfy $(T - t_i) |Rm| (x_i, t_i) \leq C$ for some $C < \infty$ independent of $i$. Hence we expect that the absolute curvatures at $(x_i, t_i)$ are not uniformly equivalent to the maximum curvatures at the times $t_i$. In particular, to obtain the desired cylindrical limits, we need to find a procedure which selects points with this property. There is such a procedure, which we discuss in the last section. In the next section, we restrict our attention to developing an intuitive understanding of how Type II singularities arise.

## 10   Degenerate neck pinch

The standard model for a Type II singularity is the degenerate neck pinch. To see how this singularity arises, imagine a family of solutions $g_s(t)$ to the Ricci flow on a topological 3-sphere parametrized by the unit interval: $s \in [0, 1]$. When $s = 0$, let the initial metric be a symmetric dumbbell with two equally sized spherical regions joined by a thin neck.

Under the Ricci flow, we expect the neck to pinch off at the center after a short time. On the other hand, when $s = 1$, let the initial metric be a lopsided dumbbell where the neck is fat and short and one of the spherical regions is smaller than the other. In this case, the smoothing effect[5] of the Ricci flow should be strong enough to make the neck not shrink and allow the smaller spherical region to pull through so that eventually the metric becomes positively curved (so that it shrinks to a point while approaching constant curvature.) Now let the initial metrics $g_s(0)$ be dumbbells where the necks become shorter and fatter and one of the spherical region becomes smaller as $s$ increases from 0 to 1. For each $s \in [0, 1]$, the solution exists up to a time $T_s < \infty$ when a singularity forms. By the continuous dependence of the solution on the initial metric, we expect that there is some parameter $s_0 \in [0, 1]$ such that for $s \in [0, s_0)$, the neck pinches off at time $T_s$, and for $s \in (s_0, 1]$, the necks do not pinch off and the metrics shrink to a point while approaching constant curvature at time $T_s$. For the solution $g_{s_0}(t)$ we expect that the neck pinches off at the same time $T_{s_0}$ the smaller spherical region shrinks to a point, leaving us with a cusplike singularity which is known as a degenerate neck pinch.

More generally, we can ask the following question. Given a compact differentiable manifold and a one-parameter family of initial metrics $g^s$, $s \in [1, 2]$ such that $g^1$ forms a Type I $S^3/\Gamma$ singularity model and $g^2$ forms a Type I $\left(S^2 \times \mathbb{R}\right)/\Gamma$ singularity model, then does there exist $s \in (1,2)$ such that $g^s$ forms a Type II singularity? The guess is that at any $s \in (1,2)$ where there is a transition in the singularity model, the singularity is Type II.

---

[5]See [BMR], [B], [S1] and [S2] for derivative estimates which hold for the Ricci flow.

The reader may be thinking to herself or himself that is all very well, but do degenerate neck pinches really exist? For the mean curvature flow of a surface in euclidean 3-space, Angenent and Velazquez [AV] have proved the existence of a degenerate neck pinch. Although this is still unknown for the Ricci flow, since the mean curvature flow is in many respects very similar to the Ricci flow, it is strongly believed that degenerate neck pinches exist for the Ricci flow. In any case, the example of the degenerate neck pinch is used to develop intuition for understanding Type II singularities. It is quite likely that its analysis is not necessary for carrying through Hamilton's program. This is not to say that the study of degenerate neck pinches is not important. Degenerate neck pinches for the Ricci flow, if they can be proved to exist and their asymptotic behavior analyzed, could be used to the help formulate conjectures concerning the analytic behavior of singularities.

For now the question is, what do we see when we take a limit of dilations about a degenerate neck pinch?

## 11 Dilations revisited

We now describe what we expect to see as the limit of dilations about any Type II singularity while keeping in mind the conjectured degenerate neck pinch as the primary example. The first thing to do is to be more careful in how we choose our sequence of points and times. As for Type I singularities, we choose $(x_i, t_i) \in M \times [0, T)$ with $t_i \nearrow T$ such that $K_i \doteq |Rm|(x_i, t_i) \geq c \cdot \max_{M \times [0, t_i]} |Rm|$ for some constant $c > 0$ independent of $i$. However, unlike in the Type I case where it is not possible, in the Type II case, we want the limit solution to satisfy

$$|Rm_\infty|(x_\infty, 0) = \sup_{M_\infty \times (-\infty, \infty)} |Rm_\infty|. \tag{1}$$

This is always possible for Type II singularities.

Let's assume we have an injectivity radius estimate so that we can obtain a complete limit solution $\left(M_\infty^3, g_\infty(t), x_\infty\right)$ with the above property. This limit solution exists for all time $t \in (-\infty, \infty)$ (such solutions are called *eternal*.) Recall that the limit solution has nonnegative sectional curvature and is not flat since $|Rm_\infty|(x_\infty, 0) = 1$. (This follows from an estimate of Hamilton [H7] and Ivey [Iv].) Hence $\left(M_\infty^3, g_\infty(t)\right)$ is either the quotient of the product of a complete solution $\left(N^2, h(t)\right)$ on a surface with $\mathbb{R}$ or has positive sectional curvature. By an application of the matrix Harnack inequality of Li-Yau type[6] [H4] and the strong maximum principle [H2] to equation (1), Hamilton has shown that any such solution is time independent (modulo the action of the diffeomorphism group on the space of metrics.) Such solutions are expected to be very special. In particular, in the first case, Hamilton has shown that $\left(N^2, h(t)\right)$ is isometric to the cigar soliton given by $\left(\mathbb{R}^2, \frac{dx^2 + dy^2}{1 + x^2 + y^2}\right)$. As we remarked earlier, this limit is not expected to occur. In the second case, which we do expect to occur, $\left(M_\infty^3, g_\infty(t)\right)$ is a complete time-independent solution with positive sectional curvature, so that $M_\infty^3 \cong \mathbb{R}^3$.

---

[6]Such differential Harnack inequalities were first proved by Li and Yau [LY] for positive solutions of the heat equation and are also known as Li-Yau or Li-Yau-Hamilton inequalities. See [ChCh] for a geometric interpretation of Hamilton's matrix inequality. See also [N] for some very recent work on linear-type inequalities.

The claim is that, if the cigar soliton never appears as a factor of a limit of dilations, then $\left(M_\infty^3, g_\infty(t)\right)$ looks like the round cylinder $S^2 \times \mathbb{R}^1$ at infinity (compare this with the remarks in section 6.) This claim is proved by what is called 'dimension reduction' where one starts with the limit solution $\left(M_\infty^3, g_\infty(t)\right)$, chooses an appropriate sequence of points $(y_i, u_i) \in M_\infty^3 \times (-\infty, \infty)$ with $y_i$ tending to spatial infinity, and obtains a new limit $\left(P_\infty^3, k_\infty(t)\right)$ of dilations about $(y_i, u_i)$. This new limit exists by the local injectivity radius estimate for odd-dimensional solitons and is shown, after possibly changing the sequence of times, to be the product of a complete solution $\left(Q^2, h(t)\right)$ on a surface with $\mathbb{R}$, where $\left(Q^2, h(t)\right)$ is either a homothetically shrinking solution or a time-independent solution. In the first case, changing the sequence of times is not necessary, and we obtain a shrinking round $S^2$ as desired. In the second case, we obtain the cigar soliton, and we hope this case does not occur. Now let's compare all of this to the example of a degenerate neck pinch.

For a degenerate neck pinch we expect that at times $t$ near the singularity time, the maximum curvature of $\left(M^3, g(t)\right)$ is attained at the tip of the degenerate neck. If we dilate about any sequence of points and times $(x_i, t_i)$ such that $t_i \nearrow T$ and $|Rm|(x_i, t_i) = \max_{M \times \{t_i\}} |Rm|$, we expect that the limit is a complete, noncompact, time-independent solution $\left(M_\infty^3, g_\infty(t)\right)$ to the Ricci flow with positive sectional curvature and which is rotationally symmetric. The existence and uniqueness of such a solution to the Ricci flow has been proved by Robert Bryant [Br], who also showed that the curvature decays like the reciprocal of the distance to the origin: $|Rm_\infty(x, t)| \sim C/d_\infty(x, O, t)$. An important observation is the following. Fix a time, say $t = 0$, and take any sequence of points $y_i \in M_\infty^3$ tending to infinity and shrink the metric $g_\infty(0)$ centered at $y_i$ so that the absolute curvature at $y_i$ becomes 1. Call these shrunken metrics $g_i$. The maximum absolute curvature of $g_i$ on $M_\infty^3$ tends to infinity as $i \to \infty$ (since $|Rm_\infty|(y_i, 0) \to 0$.) However, since the curvature of $g_\infty(0)$ decays at a rate slower than $1/d_\infty(x, O, 0)^2$, the curvature of $g_i$ stays uniformly bounded on larger and larger balls centered at $y_i$. That is, the high curvature region of $\left(M_\infty^3, g_i\right)$ tends to spatial infinity with respect to the origins $y_i$. Hence, by the compactness theorem, we can take a pointed limit of some subsequence: $\left(M_\infty^3, g_i, y_i\right) \to \left(P_\infty^3, k_\infty, y_\infty\right)$. This limit is a shrinking round product cylinder. We obtained $\left(P_\infty^3, k_\infty, y_\infty\right)$ as the limit of a limit of dilations of the original solution $\left(M^3, g(t)\right)$. This means that $\left(P_\infty^3, k_\infty, y_\infty\right)$ itself is a limit. Indeed, all we have to do is adjust the points and take a subsequence of the times and consider the sequence $(y_i, t_{j_i})$. Actually $y_i$ is a point in $M_\infty^3$, but by choosing the sequence $\{j_i\}$ of positive integers tending to infinity appropriately, not only do we guarantee that $y_i$ corresponds to a point in $M^3$, but also that $\left(M^3, g\left(t_{j_i}\right)\right)$ is closer and closer to $\left(M_\infty^3, g_\infty\left(0\right)\right)$ in larger and larger balls centered at $y_i$. This implies that $\left(P_\infty^3, k_\infty, y_\infty\right)$ is a limit of dilations of $\left(M^3, g\left(t_{j_i}\right)\right)$ about $y_i$ as we claimed.

The above argument is just the technique of dimension reduction which we mentioned in the previous paragraph. One often begins with a complete, noncompact, time-independent solution $\left(M_\infty^3, g_\infty\left(t\right)\right)$ with positive sectional curvature (but not necessarily rotationally symmetric) and shows that there exists a sequence of points $y_i$ tending to infinity such that if we shrink the metric $g_\infty\left(0\right)$ centered at $y_i$ so that the absolute curvature at $y_i$ becomes 1, then there exists a subsequence which converges to a solution to the Ricci flow which is the product of a solution $\left(N_\infty^2, h_\infty\left(t\right)\right)$ to the Ricci flow on a surface with $\mathbb{R}$. What is needed for this is to work that the curvature of $\left(M_\infty^3, g_\infty\left(0\right)\right)$ does not decay quadratically or faster: $A \doteq \limsup_{x \to \infty} |Rm_\infty\left(x, 0\right)| \, d_\infty\left(x, O, 0\right)^2 = \infty$. To obtain a limit where the surface solution is either a round shrinking sphere or the cigar soliton, one needs to take a further limit of the limit solution $\left(M_\infty^3, g_\infty\left(t\right)\right)$ where one dilates about an appropriate sequence of points $(z_i, \tau_i)$ where $\tau_i \to -\infty$. More generally, one can perform dimension reduction whenever $A = \infty$ provided there is a local injectivity radius estimate.

## 12 Little Loop Lemma

The Little Loop Lemma [H7] says:

**Theorem 1 (Hamilton)** *Let $\left(\mathcal{M}^n, g\left(t\right)\right)$, $t \in [0, T)$, be a solution to the Ricci flow on a compact Riemannian manifold with nonnegative curvature operator $Rm \geq 0$. There exists $\varepsilon > 0$ and $\delta > 0$ such that if $(x_0, t_0) \in \mathcal{M} \times [0, T)$ is a point and time satisfying*

$$|Rm\left(x, t_0\right)| \leq \frac{\varepsilon}{W^2} \quad on \ B\left(x_0, W\right)$$

*for some $W > 0$, then*

$$inj(x_0, t_0) \geq \delta W.$$

*We can choose $\varepsilon > 0$ to depend only on $n$.*

There is a gap, acknowledged by Hamilton, in the proof of this result in [H7]. In particular, the statement:

> Since $s$ is conelike at $P_*$ but $R$ is smooth, the minimum is not at $P_*$, so $s > 0$.

in the proof of Lemma 15.7 is not true. However, Hamilton and Yau have announced a rigorous proof of Theorem 1 when $n = 3$.

**Conjecture 2** *(Hamilton) In the Little Loop Lemma, if $n = 3$, we can remove the hypothesis that $Rm \geq 0$.*

See [P], sections 4 for an apparent proof of this. Now Theorem 1 is related to injectivity radius estimates for solutions to the Ricci flow.

**Definition 3** *A solution $(\mathcal{M}^n, g(t))$, $t \in [0, T)$, satisfies an **injectivity radius estimate** if there exists a constant $c > 0$ such that*

$$inj(x, t)^2 \geq \frac{c}{\max_{y \in \mathcal{M}} |Rm(y, t)|}$$

*for all $x \in \mathcal{M}$ and $t \in [0, T)$.*

An injectivity radius estimate implies that one can apply the Gromov-type compactness theorem to obtain limits of dilations about a singularity.

**Lemma 4** *If $(\mathcal{M}^n, g(t))$, $t \in [0, T)$, is a solution which satisfies the conclusion of the Little Loop Lemma, then an injectivity radius estimate holds for this solution.*

**Proof.** Suppose an injectivity radius estimate does not hold. Then there exists a sequence of points $(x_i, t_i) \in \mathcal{M} \times [0, T)$ with $t_i \nearrow T$ such that

$$\text{inj}(x_i, t_i)^2 \max_{y \in \mathcal{M}} |Rm(y, t_i)| \doteq \mu_i \to 0.$$

Define

$$W \doteq \mu_i^{-1/4} \text{inj}(x_i, t_i) = \mu_i^{1/4} \max_{y \in \mathcal{M}} |Rm(y, t_i)|^{-1/2}.$$

We have

$$\text{inj}(x_i, t_i) = \mu_i^{1/4} W$$

and

$$|Rm(x, t_i)| \leq \max_{y \in \mathcal{M}} |Rm(y, t_i)| = \frac{\mu_i^{1/2}}{W^2} \quad \text{for all } x \in \mathcal{M}.$$

Since $\mu_i^{1/4}$ and $\mu_i^{1/2}$ both tend to zero, this contradicts the conclusion of the Little Loop Lemma. ∎